\expandafter\ifx\csname mthreemacsloaded\endcsname\relax\else \fi

\magnification1100
\input amstex


 \catcode`\@=11
 \let\wlog@ld\wlog
 \def\wlog#1{\relax}

 \newif\ifIN@
 \def\m@rker{\m@@rker}
 \def\IN@{\expandafter\INN@\expandafter}
 \long\def\INN@0#1@#2@{\long\def\NI@##1#1##2##3\ENDNI@
    {\ifx\m@rker##2\IN@false\else\IN@true\fi}%
     \expandafter\NI@#2@@#1\m@rker\ENDNI@}
  \newtoks\Initialtoks@  \newtoks\Terminaltoks@
  \def\SPLIT@{\expandafter\SPLITT@\expandafter}
  \def\SPLITT@0#1@#2@{\def\TTILPS@##1#1##2@{%
     \Initialtoks@{##1}\Terminaltoks@{##2}}\expandafter\TTILPS@#2@}
  \newtoks\Trimtoks@

 \def\ForeTrim@{\expandafter\ForeTrim@@\expandafter}
 \def\ForePrim@0 #1@{\Trimtoks@{#1}}
 \def\ForeTrim@@0#1@{\IN@0\m@rker. @\m@rker.#1@%
     \ifIN@\ForePrim@0#1@%
     \else\Trimtoks@\expandafter{#1}\fi}
 
  \def\Trim@0#1@{%
      \ForeTrim@0#1@%
      \IN@0 @\the\Trimtoks@ @%
        \ifIN@
             \SPLIT@0 @\the\Trimtoks@ @\Trimtoks@\Initialtoks@
             \IN@0\the\Terminaltoks@ @ @%
                 \ifIN@
                 \else \Trimtoks@ {FigNameWithSpace}%
                 \fi
        \fi
      }

  \font\titlebold=cmbx12 scaled 1200
  \font\twelvebold=cmbx12
  \font\tenbold=cmbx10
  \font\ninebold=cmbx9
  \font\sevenbold=cmbx7
  \font\fivebold=cmbx5

  \input amssym.def \input amssym
     \font\titlemsa=msam10 at 14.4pt
     \font\titlemsb=msbm10 at 14.4pt
     \font\titleeufm=eufm10 at 14.4pt
     \font\twelvemsa=msam10 scaled 1200
     \font\twelvemsb=msbm10 scaled 1200
     \font\twelveeufm=eufm10 scaled 1200
     \font\ninemsa=msam9
     \font\ninemsb=msbm9
     \font\nineeufm=eufm9

   \ifx\cyrfam\undefined
   \else
     \immediate\write16{}%
     \message{ !!! cyr fonts already defined. !!! }
     \message{ --- edit out superfluous font defs? }
   \fi
   \newfam\cyrfam
       \font\titlecyr=wncyr10 scaled 1440 
       \font\twelvecyr=wncyr10 scaled 1200
       \font\tencyr=wncyr10
       \font\ninecyr=wncyr9
       \font\sevencyr=wncyr7
       \font\sixcyr=wncyr6

   \newfam\eusmfam
       \font\titleeusm=eusm10 scaled 1440
       \font\twelveeusm=eusm10 scaled 1200
       \font\teneusm=eusm10
       \font\nineeusm=eusm9
       \font\seveneusm=eusm7
       
       \font\fiveeusm=eusm5

\let\Cal\cal

    \font\ninemrm=cmr9 
    \font\ninei=cmmi9
    \font\ninesy=cmsy9 
    \skewchar\ninei='177
    \skewchar\ninesy='60

  \font\twelvemrm=cmr10 at 12pt 
  \font\twelvei=cmmi10 at 12pt
  \font\twelvesy=cmsy10 at 12pt

  \font\titlemrm=cmr10 at 14.4pt 
  \font\titlei=cmmi10 at 14.4pt
  \font\titlesy=cmsy10 at 14.4pt


  \def\Smallfonts{\ninepoint}

  \def\Hfont{\titlepoint\bf}
  \def\Authorfont{\twelvepoint\it}
  \def\HHfont{\twelvepoint\bf}
  \def\HHHfont{\bf}
  \def\Bibfont{\tenbf}
  \def\Coordfont{\nineit }

  \def \thfont {\bf }
  \def \pffont {\it\itSpacing }
  \def \rkfont {\bf }
  \def \dffont {\bf }
  \def \egfont {\bf }

 \def\ninepoint{%
  \def\rm{\fam0\ninerm}%
    \textfont0=\ninemrm  \scriptfont0=\sevenrm  \scriptscriptfont0=\fiverm
    \textfont1=\ninei    \scriptfont1=\seveni   \scriptscriptfont1=\fivei
  \def\mit{\fam1\ninei}%
  \def\oldstyle{\fam1\ninei}%
    \textfont2=\ninesy   \scriptfont2=\sevensy  \scriptscriptfont2=\fivesy
    \textfont3=\tenex    \scriptfont3=\tenex    \scriptscriptfont3=\tenex
  \def\it{\fam\itfam\nineit}%
    \textfont\itfam=\nineit
  \def\bf{\ifmmode\fam\bffam\else\ninebf\fi}%
    \textfont\bffam=\ninebold 
    \scriptfont\bffam=\sevenbold 
    \scriptscriptfont\bffam=\fivebold%
  \def\msa{\fam\msafam\ninemsa}%
    \textfont\msafam=\ninemsa 
    \scriptfont\msafam=\sevenmsa
    \scriptscriptfont\msafam=\fivemsa%
  \def\msb{\fam\msbfam\ninemsb}%
    \textfont\msbfam=\ninemsb%
    \scriptfont\msbfam=\sevenmsb%
    \scriptscriptfont\msbfam=\fivemsb%
  \def\eufm{\fam\eufmfam\nineeufm}%
    \textfont\eufmfam=\nineeufm
    \scriptfont\eufmfam=\seveneufm
    \scriptscriptfont\eufmfam=\fiveeufm
   \def\eusm{\fam\eusmfam\nineeusm}%
     \textfont\eusmfam=\nineeusm
     \scriptfont\eusmfam=\seveneusm
     \scriptscriptfont\eusmfam=\fiveeusm
   \def\cyr{\fam\cyrfam\ninecyr}%
     \textfont\cyrfam=\ninecyr
     \scriptfont\cyrfam=\sevencyr
     \scriptscriptfont\cyrfam=\sixcyr
  \setbox\strutbox=\hbox{\vrule
      height7pt depth3pt width0pt}%
   \baselineskip=10.8pt\rm}

 \let\eightpoint\ninepoint 

 \def\tenpoint{%
  \def\rm{\fam0\tenrm}%
    \textfont0=\tenmrm \scriptfont0=\sevenrm \scriptscriptfont0=\fiverm%
  \def\mit{\fam1\teni}%
  \def\oldstyle{\fam1\teni}%
    \textfont1=\teni   \scriptfont1=\seveni  \scriptscriptfont1=\fivei%
    \textfont2=\tensy  \scriptfont2=\sevensy \scriptscriptfont2=\fivesy%
    \textfont3=\tenex  \scriptfont3=\tenex   \scriptscriptfont3=\tenex%
  \def\it{\fam\itfam\tenit}%
    \textfont\itfam=\tenit%
  \def\bf{\ifmmode\fam\bffam\else\tenbf\fi}%
    \textfont\bffam=\tenbold
    \scriptfont\bffam=\sevenbold%
    \scriptscriptfont\bffam=\fivebold%
  \def\msa{\fam\msafam\tenmsa}%
    \textfont\msafam=\tenmsa%
    \scriptfont\msafam=\sevenmsa%
    \scriptscriptfont\msafam=\fivemsa%
  \def\msb{\fam\msbfam\tenmsb}%
    \textfont\msbfam=\tenmsb%
    \scriptfont\msbfam=\sevenmsb%
    \scriptscriptfont\msbfam=\fivemsb%
  \def\eufm{\fam\eufmfam\teneufm}%
   \textfont\eufmfam=\teneufm
   \scriptfont\eufmfam=\seveneufm
   \scriptscriptfont\eufmfam=\fiveeufm
   \def\eusm{\fam\eusmfam\teneusm}%
    \textfont\eusmfam=\teneusm
    \scriptfont\eusmfam=\seveneusm
    \scriptscriptfont\eusmfam=\fiveeusm
   \def\cyr{\fam\cyrfam\tencyr}%
    \textfont\cyrfam=\tencyr
    \scriptfont\cyrfam=\sevencyr
    \scriptscriptfont\cyrfam=\sixcyr
  \setbox\strutbox=\hbox{\vrule %
      height8.5pt depth3.5ptwidth0pt}%
  \baselineskip=\StdBaselineskip\rm}

 \def\twelvepoint{%
  \def\rm{\fam0\twelverm}%
    \textfont0=\twelvemrm \scriptfont0=\tenmrm \scriptscriptfont0=\sevenrm
    \textfont1=\twelvei   \scriptfont1=\teni   \scriptscriptfont1=\seveni
  \def\mit{\fam1\twelvei}%
  \def\oldstyle{\fam1\twelvei}%
    \textfont2=\twelvesy  \scriptfont2=\tensy  \scriptscriptfont2=\sevensy
    \textfont3=\tenex  \scriptfont3=\tenex  \scriptscriptfont3=\tenex
  \def\it{\fam\itfam\twelveit}%
    \textfont\itfam=\twelveit
  \def\bf{\ifmmode\fam\bffam\else\twelvebf\fi}%
    \textfont\bffam=\twelvebold
    \scriptfont\bffam=\tenbold%
    \scriptscriptfont\bffam=\sevenbold%
  \def\msa{\fam\msafam\twelvemsa}%
    \textfont\msafam=\twelvemsa%
    \scriptfont\msafam=\tenmsa%
    \scriptscriptfont\msafam=\sevenmsa%
  \def\msb{\fam\msbfam\twelvemsb}%
    \textfont\msbfam=\twelvemsb%
    \scriptfont\msbfam=\tenmsb%
    \scriptscriptfont\msbfam=\sevenmsb%
  \def\eufm{\fam\eufmfam\twelveeufm}%
   \textfont\eufmfam=\twelveeufm
   \scriptfont\eufmfam=\teneufm
   \scriptscriptfont\eufmfam=\seveneufm
   \def\eusm{\fam\eusmfam\twelveeusm}%
    \textfont\eusmfam=\twelveeusm
    \scriptfont\eusmfam=\teneusm
    \scriptscriptfont\eusmfam=\seveneusm
   \def\cyr{\fam\cyrfam\tencyr}%
    \textfont\cyrfam=\twelvecyr
    \scriptfont\cyrfam=\tencyr
    \scriptscriptfont\cyrfam=\sevencyr
  \setbox\strutbox=\hbox{\vrule
      height10.2pt depth4.55pt width0pt}%
  \baselineskip=14pt\rm}

 \def\titlepoint{%
    \textfont0=\titlemrm \scriptfont0=\twelvemrm \scriptscriptfont0=\tenmrm
    \textfont1=\titlei   \scriptfont1=\twelvei   \scriptscriptfont1=\teni
  \def\mit{\fam1\titlei}%
  \def\oldstyle{\fam1\titlei}%
    \textfont2=\titlesy  \scriptfont2=\twelvesy  \scriptscriptfont2=\tensy
    \textfont3=\tenex
    \scriptfont3=\tenex
    \scriptscriptfont3=\tenex
  \def\it{\fam\itfam\titleit}%
    \textfont\itfam=\titleit
  \def\bf{\ifmmode\fam\bffam\else\titlebf\fi}%
    \textfont\bffam=\titlebold
    \scriptfont\bffam=\twelvebold%
    \scriptscriptfont\bffam=\tenbold%
  \def\msa{\fam\msafam\titlemsa}%
    \textfont\msafam=\titlemsa%
    \scriptfont\msafam=\twelvemsa%
    \scriptscriptfont\msafam=\tenmsa%
  \def\msb{\fam\msbfam\titlemsb}%
    \textfont\msbfam=\titlemsb%
    \scriptfont\msbfam=\twelvemsb%
    \scriptscriptfont\msbfam=\tenmsb%
  \def\eufm{\fam\eufmfam\titleeufm}%
    \textfont\eufmfam=\titleeufm
    \scriptfont\eufmfam=\twelveeufm
    \scriptscriptfont\eufmfam=\teneufm
   \def\eusm{\fam\eusmfam\titleeusm}%
     \textfont\eusmfam=\titleeusm
     \scriptfont\eusmfam=\twelveeusm
     \scriptscriptfont\eusmfam=\teneusm
   \def\cyr{\fam\cyrfam\tencyr}%
    \textfont\cyrfam=\titlecyr
    \scriptfont\cyrfam=\twelvecyr
    \scriptscriptfont\cyrfam=\tencyr
  \setbox\strutbox=\hbox{\vrule
      height12.3pt depth5.54pt width0pt}%
  \baselineskip=16pt\rm}

\newbox\AuthorBox\newbox\TitleBox
\newbox\TFLinebox
\newbox\FLinebox
\newbox\HLinebox
\def\SetTFLinebox#1{\setbox\TFLinebox=\hbox{#1}}
\def\SetFLinebox#1{\setbox\FLinebox=\hbox{#1}}
\def\SetHLinebox#1{\setbox\HLinebox=\hbox{#1}}

 \def\SetAuthorHead#1{%
     \setbox\AuthorBox=\hbox{\ninepoint \it 
           \ignorespaces\frenchspacing#1\unskip}}
 \def\SetTitleHead#1{%
     \setbox\TitleBox=\hbox{\ninepoint \it
           \ignorespaces\frenchspacing#1\unskip}}

  \def\itSpacing{\relax}
  \def\itSpacingOff{\relax}


 \def\Hrule{\hrule width0pt height0pt}

  \newskip\ProcSkip \ProcSkip 8pt plus2pt minus2pt

 \newskip\LastSkip
 \def\SaveLastSkip{\LastSkip\lastskip}
 \def\RestoreLastSkip{\vskip-\LastSkip\vskip\LastSkip}

 \def\NoindentAfter{\everypar={\setbox0=\lastbox\everypar={}}}

 \long\def\H#1\par#2\par{\notenumber=0 \titlepagetrue%
    {
    \baselineskip=20pt
    \parindent=0pt\parskip=0pt\frenchspacing
    \leftskip=0pt plus .2\hsize minus .3\hsize
    \rightskip=0pt plus .2\hsize minus .3\hsize
 \def\\{\unskip\break}%
    \pretolerance=10000 \Hfont #1\unskip\break
     \vskip7pt\Hrule
\hfill \Authorfont #2\hfill\hfill\unskip}
    \vskip48pt plus 4pt minus 4pt
    \par\NoindentAfter\rm}

 \long\def\Hi#1\par#2\par{\notenumber=0 \titlepagetrue%
    {  \baselineskip=0pt  \parindent=0pt\parskip=0pt\frenchspacing
    \leftskip=0pt plus .2\hsize minus .3\hsize
    \rightskip=0pt plus .2\hsize minus .3\hsize
}
    \rm}


 \newdimen\PageRemainder
  \def\SetPageRemainder{
     \PageRemainder=\pagegoal
     \ifdim\PageRemainder=\maxdimen\PageRemainder=\vsize
     \else\advance\PageRemainder by -1\pagetotal\fi}

  \def\Rpt@{}\def\Rpt@@{}

  \long\def\HH#1\par{\par
  \SaveLastSkip\removelastskip\goodbreak
  \ifdim\LastSkip<30pt 
     \LastSkip 30pt
plus 3pt minus 2pt\fi
  \SetPageRemainder\advance\PageRemainder-\LastSkip
  \ifdim\PageRemainder<150pt
       \edef\Rpt@{remain = \the\PageRemainder\noexpand\\
                pagetotal=\the\pagetotal\noexpand\\
                           pagegoal=\the\pagegoal}%
          \fi
   \ifdim\PageRemainder<65pt 
       \ifdim\PageRemainder > 0pt
          \edef\Rpt@@{\noexpand\\
                      Had HH PageRemainder$<$\relax 65pt\noexpand\\
                      Hence forced break!}%
     \vskip 0pt plus .2\PageRemainder\eject 
    \fi\fi
    \vskip\LastSkip\Hrule 
    \pretolerance=10000\rightskip=0pt plus 3em
    \hangafter1 \hangindent=2.2em%
    \noindent
    \HHfont \unskip \Ednote{\Rpt@\Rpt@@}%
            \def\Rpt@{}\def\Rpt@@{}%
            \ignorespaces
            #1\par\rightskip=0pt\pretolerance=\StdPretolerance%
    \NoindentAfter
\tenpoint\rm%
     \medskip \vskip\ProcSkip}

  \long\def\HHH#1\par{\par%
  \SaveLastSkip\removelastskip\goodbreak
  \ifdim\LastSkip<\ProcSkip%
     \LastSkip\ProcSkip\fi
  \SetPageRemainder\advance\PageRemainder-\LastSkip
  \ifdim\PageRemainder<150pt
       \edef\Rpt@{remain = \the\PageRemainder\noexpand\\
                pagetotal=\the\pagetotal\noexpand\\
                           pagegoal=\the\pagegoal}%
       \fi
   \ifdim\PageRemainder<48pt  
        \ifdim\PageRemainder > 0pt
             \edef\Rpt@@{\noexpand\\
                      Had HHH PageRemainder$<$\relax48pt\noexpand\\
                      Hence forced break!}%
       \vskip 0pt plus .2\PageRemainder\eject 
      \fi\fi
   \vskip\LastSkip\par\noindent
   \HHHfont \unskip\Ednote{\Rpt@\Rpt@@}%
  \def\Rpt@{}\def\Rpt@@{}%
  \ignorespaces
   #1\unskip.\quad\rm\ignorespaces
   \ignorepars}

  \long\def\ignorepars#1\par{\def\Test{#1}%
     \ifx\Test\Empty\def\This{\ignorepars}%
        \else\def\This{\Test\par}\fi
           \This}
  \def\Empty{}

 \def\Abstract#1\par{\bgroup\Smallfonts\narrower\HHH #1\par}
 \def\endAbstract{\par\egroup}


 \def\ProcBreak{\par%
    \ifdim\lastskip<8pt%
    \removelastskip%
    \penalty-200\vskip\ProcSkip\fi}

 \def\th#1\par{\ProcBreak \noindent
   {\thfont\ignorespaces
    #1\unskip.}\it\itSpacing\kern.4em\ignorepars}

 \def\endth{\ProcBreak\rm\itSpacingOff }


 \def\pf#1\par{\ProcBreak %
    \noindent\pffont#1\unskip.\rm\itSpacingOff{\kern .7em}\ignorepars}

 \def\endpf{\medskip \ProcBreak } 

  \def\qedbox{\hbox{\vbox{
    \hrule width0.2cm height0.2pt
    \hbox to 0.2cm{\vrule height 0.2cm width 0.2pt
             \hfil\vrule height0.2cm width 0.2pt}
    \hrule width0.2cm height 0.2pt}\kern1pt}}

  \def\qed{\ifmmode\qedbox
    \else\unskip\ \hglue0mm\hfill\qedbox\ProcBreak\fi}

  \def \rk #1\par{\ProcBreak
     \noindent{\rkfont\ignorespaces #1\unskip.}%
     \rm\kern.6em\ignorepars}

  \def \endrk {\medskip\ProcBreak }

  \def \df #1\par{\ProcBreak
     \noindent{\dffont\unskip\ignorespaces #1\unskip.}%
     \rm\kern.6em\ignorepars}

  \def \eg #1\par{\ProcBreak
     \noindent\egfont\unskip\ignorespaces #1\unskip.
     \rm\kern.6em\ignorepars}

  \def \endeg {\medskip\ProcBreak }

  \newdimen\Overhang

   \def\MaxTag@#1#2#3#4#5{\setbox0=\hbox{#4\ignorespaces#2\unskip}%
     \dimen0=\wd0\advance\dimen0 by#3
     \ifdim\dimen0<#5\relax\dimen0=#5\fi
     \expandafter\edef\csname #1Hang\endcsname{\the\dimen0}}

 \def\MaxItemTag#1{\MaxTag@{Item}{#1}{.4em}{\ItemStyle}{\parindent}}%
 \def\MaxItemItemTag#1{%
        \MaxTag@{ItemItem}{#1}{.4em}{\ItemItemStyle}{\parindent}}
 \def\MaxNrTag#1{\MaxTag@{Nr}{#1}{.5em}{\NrStyle}{\parindent}}
 \def\MaxReferenceTag#1{%
        \MaxTag@{Reference}{[#1]}{.6em}{\ninerm}{\parindent}}
 \def\MaxFootTag#1{\MaxTag@{Foot}{#1}{.4em}{\ninerm}{\z@}}

  \def\SetOverhang@{\Overhang=.8\dimen0%
     \advance\Overhang by \wd0\relax
     \ifdim\Overhang>\hangindent\relax
       \advance\Overhang by .25\dimen0%
       \Ednote{Tag is pushing text.}\osumess{Tag is pushing text.}%
     \else\Overhang=\hangindent
     \fi}

   \def\Item#1{\par\noindent
      \hangafter1\hangindent=\ItemHang
      \setbox0=\hbox{\ItemStyle\ignorespaces#1\unskip}%
      \dimen0=.4em\SetOverhang@
      \rlap{\box0}\kern\Overhang\ignorespaces}

   \def\ItemItem#1{\par\noindent
      \hangafter1\hangindent=\ItemItemHang
      \setbox0=\hbox{\ItemItemStyle\ignorespaces#1\unskip}%
      \dimen0=.4em\SetOverhang@
      \advance\hangindent by \ItemHang
      \kern\ItemHang\rlap{\box0}%
      \kern\Overhang\ignorespaces}

  \def\Nr#1{\par\noindent\hangindent=\NrHang 
    \setbox0=\hbox{\NrStyle\ignorespaces#1\unskip}%
    \dimen0=.5em\SetOverhang@
    \rlap{\box0}\kern\Overhang
    \hangindent=\z@\ignorespaces}

   \newskip\Rosterskip\Rosterskip 1pt plus1pt 
   \def\Roster{\par\ifdim\lastskip<\Rosterskip\removelastskip\vskip\Rosterskip\fi
    \bgroup}
   \def\endRoster{\par\global\edef\LastSkip@{\the\lastskip}\removelastskip
       \egroup\penalty-50\LastSkip\LastSkip@\relax
       \ifdim\LastSkip<\Rosterskip\LastSkip\Rosterskip\fi
       \vskip\LastSkip}




 \def\cite#1{
    \def\nextiii@##1,##2\end@{{\frenchspacing\rm 
      \lBr\ignorespaces##1\unskip{\rm,~\ignorespaces##2}\rBr}}%
    \IN@0,@#1@%
    \ifIN@\def\next{\nextiii@#1\end@}\else
    \def\next{{\rm\lBr#1\rBr}}\fi\next}


   \def \Bib#1\par{%
       \par\removelastskip\SetPageRemainder
       \ifdim\PageRemainder < 97pt
        \ifdim\PageRemainder > 0pt
        \vfill\eject
       \fi\fi
    \ProcBreak \par\begingroup\parskip=0 pt%
    \goodbreak \vskip 15 pt plus 10 pt
    \noindent\null\hfill\Bibfont
      \ignorespaces #1\unskip\hfill\null\par 
    \frenchspacing \Smallfonts\rm
    \parskip=2.5 pt plus 1 pt minus.5pt%
    \nobreak\vskip 12pt plus 2pt minus2pt\nobreak
    \leftskip=0 pt \baselineskip=10.5pt}

 \def\ReferenceTagSlide{0em}
  \def\ReferenceTagGap{.5em}

  \def \rf#1{\par\noindent
     \hangafter1\hangindent=\ReferenceHang      
     \setbox0=\hbox{\ninerm[\ignorespaces#1\unskip]}%
     \dimen0=\ReferenceTagGap\SetOverhang@
     \rlap{\kern\ReferenceTagSlide\box0}%
     \kern\Overhang\ignorespaces}

  \def\ref#1\par#2\par#3\par#4\par{%
     \rf{#1}#2\unskip,\ #3\unskip,\
     #4\unskip.}

  \def\endBib{\par\endgroup\vskip 12pt minus 6pt }


  \long\def\Coordinates#1\endCoordinates{
 {\par\vskip4pt\def\\{\unskip, }\Coordfont\baselineskip10.5pt\noindent#1}}

 \def\pagecontents{
  \gdef\Pagetot@l{\pagetotal}
  \ifvoid\TRMargIns\else
    \rlap{\kern\hsize\kern10pt\vbox to 0pt{%
         \box\TRMargIns\vss}}\fi
  \ifvoid\topins\else\unvbox\topins\fi
   \dimen@=\dp\@cclv \unvbox\@cclv 
   \ifvoid\footins\else 
     \vskip\skip\footins
     \footnoterule
     \unvbox\footins\fi
   \ifr@ggedbottom \kern-\dimen@ \vfil \fi}


 \newcount\Ht 

 \def \Acc{\expandafter } 

 \def\swthat{\raise -1.1 ex\hbox{\sam$\widehat{}$}}
 \def\swttilde{\raise -1.2 ex\hbox{\sam$\widetilde{}$}}
 \def \overdot{{\raise .2 ex \hbox to 0pt {\hss\bf\smash{.}\hss}}}
 \def \overcircle{{\raise .1 ex \hbox to 0pt
    {\sam$\eightpoint\scriptstyle\hss\circ\hss$}}}

 \def \Mathaccent#1#2{{\sam 
  \setbox4=\hbox{$\vphantom{#2}$}
  \Ht=\ht4 
  \setbox5=\hbox{${#1}$}
  \setbox6=\hbox{${#2}$}
  \setbox7=\hbox to .5\wd6{}
  \copy7\kern .1\Ht \raise\Ht sp\hbox{\copy5}\kern-.1\Ht
  \copy7\llap{\box6}
  }}

  \def\SwtCheck #1{
        \ifmmode \check{#1}%
                \else \v {#1}%
                \fi}

 \def\barpartial {%
   \kern .17 em
    \overline {\kern -.17 em\partial\kern-.03 em}%
    \kern .03 em}

 
  \def\Overline#1{\setbox1=\hbox{\sam ${#1}$}%
      \ifdim \wd1 > 6pt
    \kern .11 em
    \overline {\kern -.11 em#1\kern-.14 em}
    \kern .14 em
  \else
    \kern .03 em
    \overline {\kern -.03 em#1\kern-.04 em}
    \kern .04 em
  \fi}

 \def\SOverline#1{\setbox1=\hbox{\sam ${#1}$}%
      \ifdim \wd1 > 7pt
    \kern .22 em
    \overline {\kern -.22 em#1\kern-.09 em}%
    \kern .09 em
  \else
    \kern .10 em
    \overline {\kern -.10 em#1\kern-.04 em}%
    \kern .04 em
  \fi}


 \def\Underline#1{\setbox1=\hbox{\sam ${#1}$}%
      \ifdim \wd1 > 6pt
    \kern .11 em
    \underline {\kern -.11 em#1\kern-.14 em}
    \kern .14 em
  \else
    \kern .03 em
    \underline {\kern -.03 em#1\kern-.04 em}
    \kern .04 em
  \fi}

 \def\SUnderline#1{\setbox1=\hbox{\sam ${#1}$}%
      \ifdim \wd1 > 7pt
    \kern .04 em
    \underline {\kern -.04 em#1\kern-.2 em}%
    \kern .2 em
  \else
    \kern .0 em
    \underline {\kern -.0 em#1\kern-.15 em}%
    \kern .15 em
  \fi}


 \def \Blackbox
   {\leavevmode\hskip .3pt \vbox
   {\hrule height 5pt\hbox{\hskip 4.5pt}}\hskip .5pt}

 \def \XX{\Blackbox\kern.5pt\Blackbox} 

  \def\.{.\kern1pt}

    \def\Hyphen{\edef\this{\the\hyphenchar\font}%
          \hyphenchar\font=-1\char\this\hyphenchar\font=\this}

 \ifx\undefined\text
  \def\text#1{\hbox{\rm #1}}\fi 



   \everymath{}  

  \def\PassMath@@{\aftergroup\AfterMath@} 

 \let\PassMath@\PassMath@@

 \def\AfterMath@{\futurelet\next\AfterMathMole@}

 \def\AfterMathMole@{
      \ifcat\next\space
          \def\this{}
      \else
      \ifcat\next\egroup %
        \def\this{\osumess{Handset mathsurround?? ---(see dollar brace)}}%
      \else
      \def\this{\AAfterMath@}
      \fi\fi
      \this}

 \def\hyphen@{-}
 \def\paren@{)}
 \def\apostr@{'}

 \def\MSC#1{\kern-.8\mathsurround#1\kern.8\mathsurround}

 \def\AAfterMath@#1{\def\Next{#1}
    \IN@0\Next @,.;:!?\relax @%
    \ifIN@\def\this{\MSC{\Next}}%
    \else
    \ifx\Next\hyphen@\def\this{\futurelet\next\AfterHyphen@}%
    \else
    \ifx\Next\paren@\def\this{#1}%
    \else 
    \ifx\Next\apostr@\def\this{#1}%
    \else \def\this{\osumess{Handset mathsurround??}%
                 #1}\fi\fi\fi\fi
    \this}

 \def\AfterHyphen@#1{\def\Next{#1}%
   \ifx\Next\hyphen@\def\this{--}\else
   \ifcat\next\space%
   \def\this{\kern-\mathsurround\kern.05em- \Next}\else
   \def\this{\kern-\mathsurround\kern.05em\Hyphen\Next}\fi\fi\this}

 \def\sam{\mathsurround=\z@\let\PassMath@\relax}  %
 \def\mas{\mathsurround=\StdMathsurround\let\PassMath@\PassMath@@}
 
 \def\Mas{\mathsurround=\StdMathsurround
                \everymath{\PassMath@}\let\PassMath@\PassMath@@}

 \def\m@th{\mathsurround=\z@\everymath{}}

 \def\m@@th{\mathsurround=\z@\everymath={}\let\m@th\relax}

\def\underbar#1{$\setbox\z@\hbox{#1}\dp\z@\z@
      \m@th \underline{\box\z@}$\relax}

\def\mathhexbox#1#2#3{\leavevmode
  \hbox{\m@@th$\m@th \mathchar"#1#2#3$}}

\def\dots{\relax\ifmmode\ldots\else$\m@th\ldots\,$\relax\fi}

\def\dotfill{\cleaders\hbox{\m@@th$\m@th \mkern1.5mu.\mkern1.5mu$}\hfill}
\def\rightarrowfill{$\m@th\mathord-\mkern-6mu%
  \cleaders\hbox{\m@@th$\mkern-2mu\mathord-\mkern-2mu$}\hfill
  \mkern-6mu\mathord\rightarrow$\relax}
\def\leftarrowfill{$\m@th\mathord\leftarrow\mkern-6mu%
  \cleaders\hbox{\m@@th$\mkern-2mu\mathord-\mkern-2mu$}\hfill
  \mkern-6mu\mathord-$\relax}

\def\downbracefill{$\m@th\braceld\leaders\vrule\hfill\braceru
  \bracelu\leaders\vrule\hfill\bracerd$\relax}
\def\upbracefill{$\m@th\bracelu\leaders\vrule\hfill\bracerd
  \braceld\leaders\vrule\hfill\braceru$\relax}

\def\angle{{\vbox{\m@@th\ialign{$\m@th\scriptstyle##$\crcr
      \not\mathrel{\mkern14mu}\crcr
      \noalign{\nointerlineskip}
      \mkern2.5mu\leaders\hrule height.34pt\hfill\mkern2.5mu\crcr}}}}

\def\big#1{{\m@@th\hbox{$\left#1\vbox to8.5\p@{}\right.\n@space$}}}
\def\Big#1{{\m@@th\hbox{$\left#1\vbox to11.5\p@{}\right.\n@space$}}}
\def\bigg#1{{\m@@th\hbox{$\left#1\vbox to14.5\p@{}\right.\n@space$}}}
\def\Bigg#1{{\m@@th\hbox{$\left#1\vbox to17.5\p@{}\right.\n@space$}}}
\def\n@space{\nulldelimiterspace\z@ \m@th}

\def\root#1\of{\setbox\rootbox\hbox{\m@@th$\m@th\scriptscriptstyle{#1}$}
  \mathpalette\r@@t}
\def\r@@t#1#2{\setbox\z@\hbox{\m@@th$\m@th#1\sqrt{#2}$\relax}
  \dimen@\ht\z@ \advance\dimen@-\dp\z@
  \mkern5mu\raise.6\dimen@\copy\rootbox \mkern-10mu \box\z@}

\def\mathph@nt#1#2{\setbox\z@\hbox{\m@@th$\m@th#1{#2}$}\finph@nt}

\def\mathsm@sh#1#2{\setbox\z@\hbox{\m@@th$\m@th#1{#2}$}\finsm@sh}

\def\@vereq#1#2{\lower.5\p@\vbox{\m@@th\baselineskip\z@skip\lineskip-.5\p@
    \ialign{$\m@th#1\hfil##\hfil$\crcr#2\crcr=\crcr}}}

\def\mathpalette#1#2{\sam\mathchoice{#1\displaystyle{#2}}%
  {#1\textstyle{#2}}{#1\scriptstyle{#2}}{#1\scriptscriptstyle{#2}}\mas}

\def\widehat#1{\setbox\z@\hbox{\sam$#1$}%
 \ifdim\wd\z@>\tw@ em\mathaccent"0\msbfam@5B{#1}%
 \else\mathaccent"0362{#1}\fi}
\def\widetilde#1{\setbox\z@\hbox{\sam$#1$}%
 \ifdim\wd\z@>\tw@ em\mathaccent"0\msbfam@5D{#1}%
 \else\mathaccent"0365{#1}\fi}

 \def\dots{\relax{}
  \ifmmode\def\thedots{\mdots@}\else\def\thedots{\tdots@}\fi %
  \thedots}

 \let\@ldeqno\eqno\let\@ldleqno\leqno
 \def\eqno{\everymath{}\@ldeqno} \def\leqno{\everymath{}\@ldleqno}

  \let\@ldeqalignno\eqalignno
  \def\eqalignno#1{\sam\@ldeqalignno{#1}\mas}
  \let\@ldeqalign\eqalign
  \def\eqalign#1{\sam\@ldeqalign{#1}\mas}

 \def\overrightarrow#1{\vbox{\m@th\ialign{##\crcr
      \rightarrowfill\crcr\noalign{\kern-\p@\nointerlineskip}
      $\hfil\displaystyle{#1}\hfil$\crcr}}}
 \def\overleftarrow#1{\vbox{\m@th\ialign{##\crcr
      \leftarrowfill\crcr\noalign{\kern-\p@\nointerlineskip}
      $\hfil\displaystyle{#1}\hfil$\crcr}}}
 \def\overbrace#1{\mathop{\vbox{\m@th\ialign{##\crcr\noalign{\kern3\p@}
      \downbracefill\crcr\noalign{\kern3\p@\nointerlineskip}
      $\hfil\displaystyle{#1}\hfil$\crcr}}}\limits}
 \def\underbrace#1{\mathop{\vtop{\m@th\ialign{##\crcr
      $\hfil\displaystyle{#1}\hfil$\crcr\noalign{\kern3\p@\nointerlineskip}
      \upbracefill\crcr\noalign{\kern3\p@}}}}\limits}

  \let\@ldmatrix\matrix
  \let\end@ldmatrix\endmatrix
  \def\matrix{\sam\@ldmatrix}
  \def\endmatrix{\end@ldmatrix\mas}
  \let\@ldgather\gather
  \let\end@ldgather\endgather
  \def\gather{\sam\@ldgather}
  \def\endgather{\end@ldgather\mas}
  \let\@ldalign\align
  \let\end@ldalign\endalign
  \def\align{\sam\@ldalign}
  \def\endalign{\end@ldalign\mas}
  \let\@ldaligned\aligned
  \let\end@ldaligned\endaligned
  \def\aligned{\sam\@ldaligned}
  \def\endaligned{\end@ldaligned\mas}
  \let\@ldtag\tag
  \def\tag{\sam\@ldtag}
   %

   \let\MinCDArrowWidth\minCDaw@




\newskip\insertskipamount\newskip\inserthardskipamount
\insertskipamount 6pt plus2pt 
\inserthardskipamount 6pt
\def\insertskip{\vskip\insertskipamount}
\newcount\SplitTest
\def\SetSplitTest{\SplitTest\insertpenalties
  \insert\topins{\floatingpenalty1}%
  \advance\SplitTest-\insertpenalties}
\def\midinsert{\par
 \SaveLastSkip\penalty-150\SetSplitTest\RestoreLastSkip
 \ifnum\SplitTest=-1
  \@midfalse\p@gefalse\else\@midtrue\fi\@ins}
\def\@ins{\par\begingroup\setbox\z@\vbox\bgroup%
  \vglue\inserthardskipamount}
\def\endinsert{\egroup 
  \if@mid \dimen@\ht\z@ \advance\dimen@\dp\z@
    \advance\dimen@\insertskipamount
    \advance\dimen@\pagetotal\advance\dimen@-\pageshrink
    \ifdim\dimen@>\pagegoal\@midfalse\p@gefalse\fi\fi
  \if@mid%
    \ifdim\lastskip<\insertskipamount\removelastskip\insertskip\fi
    \nointerlineskip\box\z@\penalty-200\insertskip
  \else%
    \SaveLastSkip
    \insert\topins{\penalty100 
    \splittopskip\z@skip
    \splitmaxdepth\maxdimen \floatingpenalty\z@
    \ifp@ge \dimen@\dp\z@
    \vbox to\vsize{\unvbox\z@\kern-\dimen@}
    \else \box\z@\nobreak\insertskip\fi}
    \RestoreLastSkip
   \fi\endgroup}


  \newcount\notenumber
  
  \def\note{\advance\notenumber by 1
    \footnote{\the\notenumber)}}

  \newbox\footbox

  \def\footnote#1{\let\@sf\empty
    \ifhmode\edef\@sf{\spacefactor\the\spacefactor}\/\fi
    \sam${}^{\fam0 #1}$\@sf\vfootnote{#1}}%

  \def\vfootnote#1{\insert\footins\bgroup
     \interlinepenalty100 \splittopskip=1pt
     \floatingpenalty=20000
     \leftskip=0pt\rightskip=0pt%
     \parindent=.3em
     \Smallfonts\rm
     \FootItem@{#1}
     \futurelet\next\fo@t}

  \def\FootItem@#1{\par\hangafter1\hangindent=\FootHang
     \setbox0=\hbox{\ignorespaces#1\unskip}%
     \dimen0=.4em\SetOverhang@
     \noindent\rlap{\box0}\kern\Overhang\ignorespaces}


  \def\fo@t{\ifcat\bgroup\noexpand\next \let\next\f@@t
    \else\let\next\f@t\fi \next}
  \def\f@@t{\bgroup\aftergroup\@foot\let\next}
  \def\f@t#1{\baselineskip=10pt\lineskip=1pt
            \lineskiplimit=0pt #1\@foot}%
  \def\@foot{
        \hbox{\vrule height0pt depth5pt width0pt}
        \egroup}
  \skip\footins=12 pt plus 0pt minus 0pt 
  \count\footins=1000 
  \dimen\footins=8in 



 \def\osumess#1{\EdSpider{\immediate\write16{Line \the\inputlineno: #1}}}%
 \def\HideEdStuff{\gdef\EdSpider##1{}}

 \font\BigSym=cmmi10 scaled \magstep 4

 \def\change{\InLMargin{\hbox{\BigSym \char63\kern10pt}}}

 \def\beginchange{\InLMargin{\hbox{\sam\twelvepoint$\heartsuit$\kern10pt}}}

 \def\endchange{\InLMargin{\hbox{\sam\twelvepoint$\spadesuit$\kern10pt}}}

 \def\InLMargin#1{\strut\vadjust{%
     \kern-\strutdepth
     \vtop to \strutdepth{%
         \baselineskip\strutdepth
         \llap{\sam$\smash{\hbox{\EdSpider{#1}}}$}\null}}}

 \def\strutdepth{\dp\strutbox}
 \def\strutheight{\ht\strutbox}

 \def\NoteInRMargin#1{\strut\vadjust{%
     \kern-1.001\strutdepth
     \vtop to \strutdepth{%
       \baselineskip\strutdepth
       \vss\rlap{\ninepoint\unskip\hskip\hsize
         \vtop to 0pt{%
           \hsize=16em\hfuzz=\hsize
           \leftskip=10pt%
           \rightskip=0pt plus 10000pt%
           \baselineskip=9.8pt\lineskip=.2pt%
           \let\\\break
           \noindent\EdSpider{#1}\vss}%
                \kern10pt}\hbox{}}
       }}

 \def\ednote#1{\NoteInRMargin{\tentt #1}}

 \def\cbar{\InLMargin{%
      \dimen0=\strutdepth\advance\dimen0 by \lineskip
      \vrule width 3pt
      height \strutheight depth \dimen0 \kern
      3pt}}

 \def\ccbar{\InLMargin{%
      \dimen0=2\strutdepth\advance\dimen0 by 2\lineskip
      \vrule width 3pt
        height 3\strutheight depth \dimen0 \kern
      3pt}}

 \newinsert\TRMargIns
 \dimen\TRMargIns=\maxdimen

  \def\Ednote#1{\insert\TRMargIns{%
       \vbox to 0pt{\hsize=140pt\hfuzz=\hsize
           \leftskip=6pt%
           \rightskip=0pt plus 10000pt%
           \baselineskip=9.8pt\lineskip=.2pt%
           \let\\\break
           \SetPageRemainder
           \vglue540pt\vglue-\PageRemainder
           \noindent\EdSpider{\tentt #1}\vss}%
       \smallskip}}

 \def\KillEdStuff{\def\ednote##1{}\def\Ednote##1{}%
      \let\change\relax\let\beginchange\relax\let\endchange\relax
       \let\cbar\relax\let\ccbar\relax}


  \topskip=12pt
  \newskip\StdBaselineskip 
  \StdBaselineskip 12pt
  \lineskip=1.1pt
  \lineskiplimit=.8pt
  \widowpenalty=10000 
  \clubpenalty=10000  
  \abovedisplayskip=6pt plus 1pt minus 1pt
  \abovedisplayshortskip=3pt plus 1.5pt
  \belowdisplayskip=6pt plus 1pt minus 1pt
  \belowdisplayshortskip=5pt plus 1pt minus 1pt
  \hfuzz=1.5pt   

  \def\StdPretolerance{100}
  \tolerance=\StdPretolerance

  \newdimen\StdMathsurround
  \StdMathsurround=1.5pt 
  \mathsurround=\StdMathsurround
  \Mas                   

   \def\prose{\relax\hbox{\kern.6\StdMathsurround}}
  
  \def\StdParskip{0pt}    
  \parskip=\StdParskip
  \parindent=0.5cm
 

  \def\Times{ptmr  } 
  \def\TimesI{ptmri  } 
  \def\TimesB{ptmb  }
  \def\TimesBI{ptmbi  }
  \def\HelveticaN{phvrrn }

  =\Times at 10bp
  =\TimesB at 10bp
  \font\tenit=\TimesI at 10bp
  =\TimesBI at 10bp

  \font\tenmrm=cmr10  


    =\Times at 9bp 
    \font\nineit=\TimesI at 9bp 
    =\TimesB at 9bp 
    =\TimesBI at 9bp 

    =\HelveticaN at 9bp 


  =\Times at 12bp
  \font\twelveit=\TimesI at 12bp
  =\TimesB at 12bp


  \font\titleit=\TimesI at 14.4bp
  =\TimesB at 14.4bp

 \SetAuthorHead{AuthorHead} 
 \SetTitleHead{TitleHead}  


  \def\lBr{\raise.125ex\hbox{[\kern.1125ex}}
  \def\rBr{\raise.125ex\hbox{\kern.1125ex]}}

 \setbox\footbox=\hbox{\Smallfonts 2)~}



  \bgroup
  \catcode`\@=11 
  \gdef\itSpacing{%
     \xspaceskip=.31em plus.1em minus.05em \sfcode `f=2001
     \itWarning@\let\itWarning@\itWarning@@}
  \gdef\itSpacingOff{%
     \xspaceskip=0pt \sfcode `f=1000
     \let\itWarning@\relax}
   \global\let\itWarning@\relax
  \gdef\itWarning@@{\errmessage{%
  Special italic spacing already in force
  (you have probably omitted an ``endth'').
  See itSpacing macro in osuPSfnt.sty
         }}
  \egroup

 \fontdimen1\titlebf=0.0pt
 \fontdimen2\titlebf=3.6135pt
 \fontdimen3\titlebf=2.8908pt
 \fontdimen4\titlebf=1.44539pt
 \fontdimen5\titlebf=6.64882pt
 \fontdimen6\titlebf=14.45398pt
 \fontdimen7\titlebf=1.60439pt

 \fontdimen1\tenbi=0.26794pt
 \fontdimen2\tenbi=2.50937pt
 \fontdimen3\tenbi=2.00749pt
 \fontdimen4\tenbi=1.00374pt
 \fontdimen5\tenbi=4.59717pt
 \fontdimen6\tenbi=10.03749pt
 \fontdimen7\tenbi=1.11415pt

 \fontdimen1\twelverm=0.0pt
 \fontdimen2\twelverm=3.01125pt
 \fontdimen3\twelverm=2.409pt
 \fontdimen4\twelverm=1.2045pt
 \fontdimen5\twelverm=5.39615pt
 \fontdimen6\twelverm=12.045pt
 \fontdimen7\twelverm=1.33699pt

 \fontdimen1\twelveit=0.27731pt
 \fontdimen2\twelveit=3.01125pt
 \fontdimen3\twelveit=2.409pt
 \fontdimen4\twelveit=1.2045pt
 \fontdimen5\twelveit=5.37207pt
 \fontdimen6\twelveit=12.045pt
 \fontdimen7\twelveit=1.33699pt

 \fontdimen1\twelvebf=0.0pt
 \fontdimen2\twelvebf=3.01125pt
 \fontdimen3\twelvebf=2.409pt
 \fontdimen4\twelvebf=1.2045pt
 \fontdimen5\twelvebf=5.5407pt
 \fontdimen6\twelvebf=12.045pt
 \fontdimen7\twelvebf=1.33699pt

 \fontdimen1\tenrm=0.0pt
 \fontdimen2\tenrm=2.50937pt
 \fontdimen3\tenrm=2.00749pt
 \fontdimen4\tenrm=1.00374pt
 \fontdimen5\tenrm=4.49678pt
 \fontdimen6\tenrm=10.03749pt
 \fontdimen7\tenrm=1.11415pt

 \fontdimen1\tenit=0.27731pt
 \fontdimen2\tenit=2.50937pt
 \fontdimen3\tenit=2.00749pt
 \fontdimen4\tenit=1.00374pt
 \fontdimen5\tenit=4.47672pt
 \fontdimen6\tenit=10.03749pt
 \fontdimen7\tenit=1.11415pt

 \fontdimen1\tenbf=0.0pt
 \fontdimen2\tenbf=2.50937pt
 \fontdimen3\tenbf=2.00749pt
 \fontdimen4\tenbf=1.00374pt
 \fontdimen5\tenbf=4.61723pt
 \fontdimen6\tenbf=10.03749pt
 \fontdimen7\tenbf=1.11415pt

 \fontdimen1\ninerm=0.0pt
 \fontdimen2\ninerm=2.25842pt
 \fontdimen3\ninerm=1.80673pt
 \fontdimen4\ninerm=0.90337pt
 \fontdimen5\ninerm=4.0471pt
 \fontdimen6\ninerm=9.03374pt
 \fontdimen7\ninerm=1.00273pt

 \fontdimen1\nineit=0.27731pt
 \fontdimen2\nineit=2.25842pt
 \fontdimen3\nineit=1.80673pt
 \fontdimen4\nineit=0.90337pt
 \fontdimen5\nineit=4.02904pt
 \fontdimen6\nineit=9.03374pt
 \fontdimen7\nineit=1.00273pt

 \fontdimen1\ninebf=0.0pt
 \fontdimen2\ninebf=2.25842pt
 \fontdimen3\ninebf=1.80673pt
 \fontdimen4\ninebf=0.90337pt
 \fontdimen5\ninebf=4.15552pt
 \fontdimen6\ninebf=9.03374pt
 \fontdimen7\ninebf=1.00273pt


 \newcount\MaxSpaceFactor
 \MaxSpaceFactor=3000 

 \def\ItemStyle{\rm}
 \def\NrStyle{\rm}
 \def\ItemItemStyle{\rm}

 \MaxItemTag{(iii)}
 \MaxItemItemTag{(iii)}
 \MaxNrTag{(2)}
 \MaxFootTag{2)}
 \def\ReferenceHang{30pt}

 \catcode`\@=\active


\loadbold

=\Times  
=\Times scaled750
=\Times scaled650
\font\rms=\Times scaled 920 

=\TimesBI scaled 860
=\TimesI scaled 860

\textfont0=\rrm  
\scriptfont0=\erm 
\scriptscriptfont0=\srm

\def\Augment#1#2{%
    \toks0\expandafter{#1}\toks2{#2}%
    \edef#1{\the\toks0\the\toks2}}

 \font\twelverma=\Times  scaled 1200
 \font\tenrma=\Times  scaled 1000
 \font\ninerma=\Times scaled 920
 =\Times scaled 840
 \font\sevenrma=\Times scaled 760
 =\Times scaled 680
 \font\fiverma=\Times scaled 600

 \Augment\tenpoint{%
  \textfont0=\tenrma  \scriptfont0=\sevenrma  
  \scriptscriptfont0=\fiverma  }

 \Augment\ninepoint{%
  \textfont0=\ninerma  \scriptfont0=\sevenrma 
  \scriptscriptfont0=\fiverma}

 \Augment\twelvepoint{%
  \textfont0=\twelverma  \scriptfont0=\ninerma  
  \scriptscriptfont0=\sevenrma}

\mathsurround=1pt
\hsize=13.45truecm
\vsize=19.5truecm
\hoffset=1.25truecm
\voffset=2truecm
\advance\baselineskip by 2pt

\predefine\til{\~}
\def\~#1{\relax\ifmmode\widetilde{#1}\else\til{#1}\fi}

\redefine \le{\leqslant}
\redefine \ge{\geqslant}
\define \wt#1{\mathaccent"0365{#1}}
\define \wh#1{\mathaccent"0362{#1}}

\define \iss{\,\Mathaccent{\raise -.8 ex\hbox{$\widetilde{}$\kern.1em}}\rightarrow\,}

\define\Car{\mathop{\fam0 C}}

\define \id{\operatorname{\fam0 id\,}}

\define \tpp{\mathop{\fam0 top}}

\define \coker{\mathop{\fam0 coker}}
\define \kr{\mathop{\fam0 ker}}

\define \gr{\operatorname{\fam0 gr}\!}

\Mas
\HideEdStuff
\rm 
 

\def\issn{{\nineit ISSN 1464-8997 (on line) 1464-8989 (printed)}}

\def\gtp{{\nineit Published 10 December 2000: \ \copyright\ Geometry \& 
Topology Publications}}

\def\gtv3{{\nineit Geometry \& Topology Monographs, Volume 3 (2000) --
Invitation to higher local fields}}


\def\lione
{{\rms Geometry \& Topology Monographs}}

\def \litwo{{\rms Volume 3: Invitation to higher local fields
}} 

\def\tinfo #1.#2.#3-#4
{{
\noindent  {\lione} \hfill 
\par 
\vskip-1.5pt
\noindent {\litwo} \hfill
\par 
\vskip-1,5pt
\noindent {\rms Part #1, section #2, pages #3--#4} \hfill
\vskip24pt 
}}

\def\tinfos #1.#2.#3-#4
{{
\noindent  {\lione} \hfill 
\par 
\vskip-1.5pt
\noindent {\litwo} \hfill
\par 
\vskip-1.5pt
\noindent {\rms Pages #3--#4} \hfill
\vskip24pt 
}}

\def\tinfoi #1
{{
\noindent  {\lione} \hfill 
\par 
\vskip-1.5pt
\noindent {\litwo} \hfill
\par 
\vskip-1.5pt
\noindent {\rms Pages iii--xi: Introduction and contents} \hfill
\vskip26pt 
}}


  \def\titlepagehead{\hfil}

  \newif\iftitlepage\titlepagefalse
  \newif\ifblankpage\blankpagefalse
  \def\makeheadline{
     \ifblankpage{}\else%
     \iftitlepage
\vbox{\line{\vbox to 8.5pt{}
\ninerm
\copy\HLinebox \hfill
\hglue5mm\ninebf\folio 
\titlepagehead}}%
      \else
\vbox{\ifodd\pageno\rightheadline\else\leftheadline\fi}%
      \fi\vskip 12pt\fi}%
     \def\rightheadline{\line{\vbox to 8.5pt{}%
      \ninerm
\copy\TitleBox \hfill
\hglue5mm\ninebf\folio}}%
     \def\leftheadline{\line{\vbox to 8.5pt{}%
        \unskip\ninerm\unskip\ninebf\folio\hglue5mm
 \hfill \copy\AuthorBox
}}

 \footline={\ifblankpage{}\else
\iftitlepage\ninepoint\sam\hfill
\line{\vbox to 8.5pt{}
\copy\TFLinebox
\hfill
\hglue5mm 
}
            \else
\ninepoint\sam\hfill
\line{\vbox to 8.5pt{}
\copy\FLinebox
\hfill 
\hglue5mm
}
\hfil\fi\global\titlepagefalse\fi}

\def\blankpage{{\blankpagetrue\noindent\hbox to 10pt{\hss}\vfill
\pagebreak}}

\tenpoint\rm 
 

\pageno=123

\font\scr=rsfs10
\define \scS{\text{\scr S}\,}
\define \fil{\text{\rm fil}}

\tinfo I.15.123-135

\SetTFLinebox{\gtp }
\SetFLinebox{\gtv3 }
\SetHLinebox{\issn}

\H 15. On the structure of the Milnor $K$-groups\\ of complete discrete valuation fields

Jinya Nakamura

\SetAuthorHead{J. Nakamura}
\SetTitleHead{Part I. Section 15. On the structure of the Milnor $K$-groups of cdv fields
\qquad\qquad}

\HH 15.0. Introduction

For a discrete valuation field $K$ the unit group $K^{*}$ of $K$ has a 
natural decreasing filtration with respect to the valuation, 
and the graded quotients of this filtration are written in terms of the residue field. 
The Milnor $K$-group $K_q(K)$ is a generalization of the unit group 
and it also has a natural decreasing filtration defined in section 4. 
However, if $K$ is of mixed characteristic and has absolute ramification index 
greater than one, the graded quotients of this filtration are known 
in some special cases only. 

Let $K$ be a  complete discrete valuation field with  residue field $k=k_K$;  
we keep the notations of section 4. Put $v_p=v_{\Bbb Q_p}$. 

A description of $\gr_nK_q(K)$ is known in the following cases:  
\Roster
\Item{(i)}
	(Bass and Tate \cite{BT})  
	$\gr_0K_q(K) \simeq K_q(k) \oplus K_{q-1}(k)$. 
\Item{(ii)}
	(Graham \cite{G}) 
	If the characteristic of $K$ and $k$ is zero, then
	$\gr_n K_q(K) \simeq \Omega^{q-1}_k$ for all $n \ge 1$. 
\Item{(iii)}	(Bloch \cite{B}, Kato \cite{Kt1}) 
	If the characteristic of $K$  and of $k$ is $p >0$ 
then $$
		\gr_nK_q(K) \simeq \coker 
			\left( 
			\Omega_k^{q-2} \longrightarrow 
			\Omega_k^{q-1}/B_s^{q-1} \oplus \Omega_k^{q-2}/B_s^{q-2}\right)$$ 
where  $\omega \quad \longmapsto (\Car^{-s}(d\omega), (-1)^q m\Car^{-s}(\omega)$ and 	where $n \ge 1$, $s=v_p(n)$ and $m=n/p^s$.
\Item{(iv)} 
	(Bloch--Kato \cite{BK})	If $K$ is of mixed characteristic $(0,p)$, then
$$
		\gr_nK_q(K) \simeq \coker\left( 
			\Omega_k^{q-2} \longrightarrow 
			\Omega_k^{q-1}/B_s^{q-1} \oplus \Omega_k^{q-2}/B_s^{q-2}\right)$$ 
where $	\omega \quad \longmapsto (\Car^{-s}(d\omega), (-1)^q m\Car^{-s}(\omega))$ and where $1 \le  n <  ep/(p-1)$ for $e=v_K(p)$, $s=v_p(n)$ and $m=n/p^s$; 
and 
$$
\split
	&\gr_{\frac{ep}{p-1}}K_q(K) \\
	&\quad \simeq  
		\coker \left(\Omega_k^{q-2} \longrightarrow 
			\Omega_k^{q-1}/(1+a\Car)B_s^{q-1} \oplus \Omega_k^{q-2}/(1+a\Car)B_s^{q-2} \right)
\endsplit
$$ 
where
$\omega \quad \longmapsto ((1+a\Car)\Car^{-s}(d\omega), (-1)^q m(1+a\Car)\Car^{-s}(\omega))$
and where $a$ is the residue class of $p/\pi^e$ for fixed prime element of $K$, $s=v_p(ep/(p-1))$ and $m=ep/(p-1)p^s$.
 
\Item{(v)}	(Kurihara \cite{Ku1}, see also section 13) 
	If $K$ is of mixed characteristic $(0,p)$ and absolutely unramified (i.e., $v_K(p)=1$),  
	then 
	$\gr_nK_q(K) \simeq \Omega_k^{q-1}/B_{n-1}^{q-1}$ for $n \ge 1$.

\Item{(vi)}
	(Nakamura \cite{N2}) 
	If $K$ is of mixed characteristic $(0,p)$ with $p>2$ and $p \nmid e=v_K(p)$,  
	then
	$$
		\gr_nK_q(K) \simeq 
		\cases 
			\text{ as in (iv)} &\quad (1 \le  n \le  ep/(p-1)) \\
			\Omega_k^{q-1}/B_{l_n+s_n}^{q-1} &\quad (n >ep/(p-1)) 
		\endcases 
	$$
	where $l_n$ is the maximal integer which satisfies $n-l_n e \ge e/(p-1)$
	and $s_n=v_p(n-l_n e)$. 
\Item{(vii)} 
	(Kurihara \cite{Ku3}) 
	If $K_0$ is 
the fraction field of the completion of the localization $\Bbb Z_p[T]
_{(p)}$ and $K=K_0(\root{p}\of{pT})$ for a prime $p \not= 2$, then 
	$$
		\gr_nK_2(K) \simeq 
		\cases
			\text{ as in (iv)} &\quad (1 \le  n \le  p) \\
			k/k^p &\quad (n=2p) \\
			k^{p^{l-2}} &\quad (n=lp, l \ge 3) \\
			0 &\quad (\text{otherwise}).
		\endcases 
	$$
\Item{(viii)} 
	(Nakamura \cite{N1}) 
	Let $K_0$ be an absolutely unramified complete discrete valuation field of mixed characteristic $(0,p)$ with $p>2$.
	If $K=K_0(\zeta_p)(\root{p}\of{\pi})$ where $\pi$ is a prime element of $K_0(\zeta_p)$ such that
	$d\pi^{p-1}=0$ in $\Omega_{{\Cal O}_{K_0(\zeta_p)}}^1$, 
	then $\gr_nK_q(K)$ are determined for all $n \ge 1$. 
	This is complicated, so we omit the details. 
\Item{(ix)} 
	(Kahn \cite{Kh}) 
	Quotients of the  Milnor $K$-groups 
of a complete discrete valuation field $K$ with
perfect residue field are computed using symbols.
\endRoster

Recall that the group of units $U_{1,K}$ can be described
as a topological $\Bbb Z_p$-module.
As a generalization of this classical result, there is an appraoch
different from (i)-(ix) for higher local fields $K$
which uses topological convergence
and $$K_q^{\tpp}(K)=K_q(K)/\cap_{l\ge1} lK_q(K)$$ (see section~6).
It provides not only the description of $\gr_nK_q(K)$ but
 of the whole $K_q^{\tpp}(K)$ 
in characteristic $p$ (Parshin \cite{P}) and in 
characteristic 0 (Fesenko \cite{F}). A complete description of the structure
of $K_q^{\tpp}(K)$ of some higher local fields with small ramification
is given by Zhukov \cite{Z}.

\smallskip


Below we discuss (vi).

\HH 15.1. Syntomic complex and Kurihara's exponential homomorphism

\HHH 15.1.1. Syntomic complex

Let $A={\Cal O}_K$ and let $A_0$ be the subring of $A$ such that $A_0$ is a  complete discrete valuation ring 
with respect to the restriction of the valuation of $K$,  
the residue field of $A_0$  coincides with $k=k_K$ and $A_0$ is absolutely unramified. 
Let $\pi$ be a fixed prime of $K$. 
Let $B=A_0[[X]]$. 
Define 
$$
\aligned 
	{\Cal J}&=\kr [ B @>{X \mapsto \pi}>> A] \\ 
	{\Cal I}&=\kr [ B @>{X \mapsto \pi}>> A
	@>\text{mod $p$}>> A/p]={\Cal J}+pB. 
\endaligned 
$$ 
Let $D$ and $J \subset D$ be the PD-envelope and the PD-ideal with respect to $B \rightarrow A$, respectively.
Let $I \subset D$ be the PD-ideal with respect to $B \rightarrow A/p$.
Namely, 
$$ 
	D=B \left[ \frac{x^j}{j!} \,\, ; \,\,  j\ge 0, x \in {\Cal J} \right], \quad 
	J=\kr(D \rightarrow A), \quad 
	I=\kr(D \rightarrow A/p).
$$ 
Let $J^{[r]}$ (resp. $I^{[r]}$) be the $r$-th divided power, which is
the ideal of $D$ generated by
$$
	\left\{ \frac{x^j}{j!} \,\,;\,\, j\ge r, x \in {\Cal J} \right\} , ~
	\left( \text{resp.~} \left\{ \frac{x^i}{i!}\frac{p^j}{j!}
		\,\,; \,\, i+j\ge r, x \in {\Cal I} \right\} \right).
$$
Notice that $I^{[0]}=J^{[0]}=D$. 
Let $I^{[n]}=J^{[n]}=D$ for a negative $n$.
We define the complexes ${\Bbb J}^{[q]}$ and ${\Bbb I}^{[q]}$ as 
$$
\aligned
{\Bbb J}^{[q]}&=[J^{[q]} @>d>> J^{[q-1]}  \otimes_B \widehat{\Omega}_B^1
	@>d>> J^{[q-2]}\otimes_B \widehat{\Omega}_B^2 \longrightarrow \cdots ] \\
{\Bbb I}^{[q]}&=[I^{[q]} @>d>> I^{[q-1]}  \otimes_B \widehat{\Omega}_B^1
	@>d>> I^{[q-2]} \otimes_B \widehat{\Omega}_B^2 \longrightarrow \cdots ]
\endaligned
$$
where $\widehat{\Omega}_B^q$ is the $p$-adic completion of $\Omega_B^q$. 
We define ${\Bbb D}={\Bbb I}^{[0]}={\Bbb J}^{[0]}$. 
\par

Let $\Bbb T$ be a fixed set of elements of $A_0^{*}$ 
such that the residue classes of all $T \in {\Bbb T}$ in $k$ forms a $p$-base of $k$.
Let $f$ be the Frobenius endomorphism of $A_0$ such that
$f(T)=T^p$ for any $T \in {\Bbb T}$ and $f(x) \equiv x^p \mod{p}$ for any $x \in A_0$. 
We extend $f$ to $B$ by $f(X)=X^p$, and to $D$ naturally.
For $0 \le r <p$ and $0 \le s$, we get
$$
	f(J^{[r]}) \subset p^rD, \quad
	f(\widehat{\Omega}_B^s) \subset p^s \widehat{\Omega}_B^s, 
$$
since
$$
\split
	f(x^{[r]})=(x^p+py)^{[r]}=(p! x^{[p]}+py)^{[r]}=p^{[r]}((p-1)!x^{[p]}+y)^r, \\
	f\bigl( z \frac{dT_1}{T_1}\wedge \dots \wedge \frac{dT_s}{T_s} \bigr)
		=z \frac{dT_1^p}{T_1^p}\wedge \dots \wedge \frac{dT_s^p}{T_s^p}
		=z p^s \frac{dT_1}{T_1}\wedge \dots \wedge \frac{dT_s}{T_s},
\endsplit
$$
where $x \in {\Cal J}$, $y$ is an element which satisfies $f(x)=x^p+py$, 
and $T_1,\dots,T_s \in {\Bbb T}\cup \{X\}$.
Thus we can define
$$
	f_q=\frac{f}{p^q} \: J^{[r]} \otimes \widehat{\Omega}_B^{q-r} \longrightarrow 
		D \otimes \widehat{\Omega}_B^{q-r}
$$
for $0 \le r<p$.
Let $\scS(q)$ and $\scS'(q)$ be the mapping fiber complexes (cf.\ Appendix) of
$$
	{\Bbb J}^{[q]} @>{1-f_q}>>  {\Bbb D} \quad\text{and}\quad 
	{\Bbb I}^{[q]} @>{1-f_q}>> {\Bbb D} 
$$ 
respectively, for $q<p$.
For simplicity, from now to the end, we assume $p$ is large enough to treat $\scS(q)$ and $\scS'(q)$.
${\scS}(q)$ is called the {\it syntomic complex} of $A$ with respect to $B$,
 and ${\scS}'(q)$ is also called the {\it syntomic complex} of $A/p$ with respect to $B$ (cf.\ \cite{Kt2}).

\th Theorem 1 {{\rm (Kurihara \cite{Ku2})}}

There exists a subgroup $S^q$ of $H^q({\scS}(q))$
such that $U_X H^q({\scS}(q)) \simeq U_1 \widehat{K_q}(A)$
where $\widehat{K_q}(A)=\varprojlim K_q(A)/p^n$
is the $p$-adic completion of $K_q(A)$ {{\rm(}}see subsection 9.1{{\rm)}}. 
\endth

\pf Outline of the proof

Let $U_X(D\otimes\widehat{\Omega}_B^{q-1})$ be the subgroup of
$D\otimes\widehat{\Omega}_B^{q-1}$ generated by
$XD\otimes\widehat{\Omega}_B^{q-1}$, $D\otimes\widehat{\Omega}_B^{q-2}\wedge dX$ and $I\otimes\widehat{\Omega}_B^{q-1}$,
and let
$$
	S^q=U_X(D\otimes\widehat{\Omega}_B^{q-1})/
	((dD\otimes \widehat{\Omega}_B^{q-2} + (1-f_q)J\otimes \widehat{\Omega}_B^{q-1}) \cap U_X(D\otimes\widehat{\Omega}_B^{q-1})).
$$
The infinite sum $\sum_{n\ge 0} f_q^n (dx)$ converges in $D\otimes\widehat{\Omega}_B^q$
for $x \in U_X(D\otimes\widehat{\Omega}_B^{q-1})$.
Thus we get a map
$$
\split
	U_X(D\otimes\widehat{\Omega}_B^{q-1}) &\longrightarrow H^q(\scS(q)) \\
	x &\longmapsto \bigl(x,\sum_{n=0}^{\infty} f_q^n (dx) \bigr)
\endsplit
$$
and we may assume $S^q$ is a subgroup of $H^q(\scS(q))$.
Let $E_q$ be the map
$$
\split
	E_q \: U_X(D\otimes \widehat{\Omega}_B^{q-1}) &\longrightarrow \widehat{K}_q(A) \\
	x\frac{dT_1}{T_1}\wedge \dots \wedge \frac{dT_{q-1}}{T_{q-1}} &\longmapsto \{E_1(x),T_1,\dots,T_{q-1}\},
\endsplit
$$
where $E_1(x)=\exp\circ(\sum_{n\ge 0}f_1^n)(x)$ is {\it Artin--Hasse's exponential homomorphism}.
In \cite{Ku2} it was shown that $E_q$ vanishes on
$$(dD\otimes \widehat{\Omega}_B^{q-2} + (1-f_q)J\otimes \widehat{\Omega}_B^{q-1}) \cap U_X(D\otimes\widehat{\Omega}_B^{q-1}),$$
hence we get the map 
$$
	E_q \: S^q \longrightarrow \widehat{K}_q(A).
$$
The image of $E_q$ coincides with $U_1\widehat{K}_q(A)$ by definition.
\par

On the other hand, define  $s_q \: \widehat{K}_q(A) \longrightarrow S^q$
by 
$$
\aligned
	& s_q(\{a_1,\dots,a_q\}) \\
		&=\sum_{i=1}^q (-1)^{i-1}\frac{1}{p}
		\log\bigl(\!\frac{f(\widetilde{a_i})}
{\widetilde{a_i}^p}\!\bigr)
		\frac{d\widetilde{a_1}}{\widetilde{a_1}} \wedge \dots \wedge
		\frac{d\widetilde{a_{i-1}}}{\widetilde{a_{i-1}}} \wedge
	f_1\bigl(\!\frac{d\widetilde{a_{i+1}}}{\widetilde{a_{i+1}}}\!\bigr) \wedge \dots \wedge
		f_1\bigl(\!\frac{d\widetilde{a_{q}}}{\widetilde{a_{q}}}\!\bigr)
\endaligned 
$$
(cf.\ \cite{Kt2}, compare with the series $\Phi$ in subsection 8.3),
where $\widetilde{a}$ is a lifting of $a$ to $D$.
One can check that $s_q \circ E_q =-\id$.
Hence $S^q \simeq U_1\widehat{K}_q(A)$.
Note that if $\zeta_p \in K$, then one can show $U_1\widehat{K}_q(A) \simeq U_1\widehat{K}_q(K)$
(see \cite{Ku4} or \cite{N2}), thus we have $S^q \simeq U_1 \widehat{K}_q(K)$.
\qed\endpf 

\eg Example

We shall prove the equality $s_q \circ E_q =-\id$ in the following simple case.
Let $q=2$.
Take an element $adT/T \in U_X(D \otimes \widehat{\Omega}_B^{q-1})$ for $T \in {\Bbb T}\cup \{X\}$.
Then 
$$
\aligned
	&s_q\circ E_q \bigl(a\frac{dT}{T}\bigr) \\
	&=s_q(\{E_1(\widetilde{a}),T\}) \\
	&=\frac{1}{p}\log\biggl( \frac{f(E_1(a))}{E_1(a)^p} \biggr) f_1\biggl(\frac{dT}{T}\biggr) \\
	&=\frac{1}{p}\biggl( \log\circ f\circ \exp \circ \sum_{n\ge 0} f_1^n(a) 
		-p\log \circ \exp \circ \sum_{n\ge 0} f_1^n(a) \biggr) \frac{dT}{T} \\
	&=\biggl( f_1 \sum_{n\ge 0} f_1^n(a) 
		-\sum_{n\ge 0} f_1^n(a) \biggr) \frac{dT}{T} \\
	&=-a \frac{dT}{T}.
\endaligned 
$$
\endeg

\HHH 15.1.2. Exponential Homomorphism

The usual exponential homomorphism
$$
\split
	\exp_{\eta} \: A &\longrightarrow A^* \\
	x &\longmapsto \exp(\eta x)=\sum_{n\ge 0} \frac{x^n}{n!}
\endsplit
$$
is defined for $\eta \in A$ such that $v_A(\eta)>e/(p-1)$.  
This map is injective.
Section 9 contains a definition of the map
$$
\aligned 
	\exp_{\eta} \: \widehat{\Omega}_A^{q-1} &\longrightarrow \widehat{K}_q(A) \\
	x\frac{dy_1}{y_1}\wedge \dots \wedge \frac{dy_{q-1}}{y_{q-1}} 
		&\longmapsto \{\exp(\eta x),y_1,\dots,y_{q-1}\}
\endaligned 
$$
for $\eta \in A$ such that $v_A(\eta)\ge 2e/(p-1)$. This map is not injective in general.
Here is a description of the kernel of $\exp_{\eta}$.

\th Theorem 2

The following sequence is exact:
$$
	H^{q-1}(\scS'(q)) \overset{\psi}\to{\longrightarrow} \Omega_A^{q-1}/pd\widehat{\Omega}_A^{q-2}
		\overset{\exp_p}\to{\longrightarrow} \widehat{K}_q(A). \tag{*}
$$
\endth

\pf  Sketch of the proof

There is an exact sequence of complexes 
$$
\alignat{3}
	0 \rightarrow 
		&\text{MF} 
		\left(\quad\,\,\CD{\Bbb J}^{[q]} \\ @V 1-f_q VV \\ {\Bbb D} \endCD \,\,\right)
		\rightarrow 
		&&\text{MF} 
	 \left(\quad\,\,\CD {\Bbb I}^{[q]} \\ @V1-f_q VV \\ {\Bbb D} \endCD\,\, \right) 
		\rightarrow
		&& {\Bbb I}^{[q]}/{\Bbb J}^{[q]}
		\rightarrow
	0, \\
	& \qquad \quad \Vert && \qquad \quad \Vert \\
	& \qquad \scS(q) && \qquad \scS'(q) 
\endalignat
$$
where MF means the mapping fiber complex.
Thus, taking cohomologies we have the following diagram with the exact top row
$$
\CD
	H^{q-1}({\scS}'(q)) @>\psi>> H^{q-1}({\Bbb I}^{[q]}/{\Bbb J}^{[q]}) @>\delta >> H^q({\scS}(q)) \\
	@. @A\text{(1)}AA @A\text{Thm.1}AA \\ 
	@. \widehat{\Omega}_A^{q-1}/pd\widehat{\Omega}_A^{q-2} @>\exp_p>> U_1 \widehat{K}_q(A),
\endCD
$$
where the map (1) is induced by
$$
	\widehat{\Omega}_A^{q-1}
	\ni \omega \longmapsto p\widetilde{\omega} \in I\otimes\widehat{\Omega}_B^{q-1}/J\otimes\widehat{\Omega}_B^{q-1}
		=({\Bbb I}^{[q]}/{\Bbb J}^{[q]})^{q-1}.
$$
We denoted the left horizontal arrow of the top row by $\psi$ and 
the right horizontal arrow of the top row by $\delta$.
The right vertical arrow is injective, thus the claims are
\Roster
	\Item{(1)} is an isomorphism,
	\Item{(2)} this diagram is commutative.
\endRoster

First we shall show (1).
Recall that 
$$
	H^{q-1}({\Bbb I}^{[q]}/{\Bbb J}^{[q]})
	=\coker \left( \frac{I^{[2]}\otimes\widehat{\Omega}_B^{q-2}}{J^{[2]}\otimes\widehat{\Omega}_B^{q-2}}
		\longrightarrow \frac{I\otimes\widehat{\Omega}_B^{q-2}}{J\otimes\widehat{\Omega}_B^{q-2}} \right).
$$
From the exact sequence
$$
	0 \longrightarrow J \longrightarrow D \longrightarrow A \longrightarrow 0,
$$
we get $D\otimes \widehat{\Omega}_B^{q-1}/J\otimes \widehat{\Omega}_B^{q-1}=A \otimes \widehat{\Omega}_B^{q-1}$
and its subgroup $I\otimes\widehat{\Omega}_B^{q-2}/J\otimes\widehat{\Omega}_B^{q-2}$ is $pA \otimes \widehat{\Omega}_B^{q-1}$
in $A\otimes \widehat{\Omega}_B^{q-1}$.
The image of $I^{[2]}\otimes \widehat{\Omega}_B^{q-2}$ in $pA \otimes \widehat{\Omega}_B^{q-1}$ is
equal to the image of 
$$
	{\Cal I}^2\otimes \widehat{\Omega}_B^{q-2}={\Cal J}^2\otimes \widehat{\Omega}_B^{q-2}
	+p{\Cal J}\widehat{\Omega}_B^{q-2} +p^2\widehat{\Omega}_B^{q-2}.
$$
On the other hand, from the exact sequence
$$
	0 \longrightarrow {\Cal J} \longrightarrow B \longrightarrow A \longrightarrow 0,
$$
we get an exact sequence
$$
	({\Cal J}/{\Cal J}^2)\otimes \widehat{\Omega}_B^{q-2}
	\overset{d}\to{\longrightarrow} A\otimes \widehat{\Omega}_B^{q-1} \longrightarrow \widehat{\Omega}_A^{q-1} \longrightarrow 0.
$$
Thus $d{\Cal J}^2\otimes \widehat{\Omega}_B^{q-2}$  vanishes on $pA\otimes \widehat{\Omega}_B^{q-1}$,
hence
$$
	H^{q-1}({\Bbb I}^{[q]}/{\Bbb J}^{[q]})
	=\frac{pA\otimes\widehat{\Omega}_B^{q-1}}{pd{\Cal J}\widehat{\Omega}_B^{q-2} +p^2d\widehat{\Omega}_B^{q-2}}
	\overset{\,p^{-1}}\to{\simeq}\frac{A\otimes\widehat{\Omega}_B^{q-1}}{d{\Cal J}\widehat{\Omega}_B^{q-2} +pd\widehat{\Omega}_B^{q-2}}
	\simeq \widehat{\Omega}_A^{q-1}/pd\widehat{\Omega}_A^{q-2},
$$
which completes the proof of (1).

Next, we shall demonstrate the commutativity of the diagram on a simple example.
Consider the case where $q=2$ and take $adT/T \in \widehat{\Omega}_A^1$ for $T \in {\Bbb T}\cup \{\pi\}$.
We want to show that the composite of
$$
	\widehat{\Omega}_A^1/pdA \overset{\text{(1)}}\to{\longrightarrow} 
	H^1({\Bbb I}^{[2]}/{\Bbb J}^{[2]}) \overset{\delta}\to{\longrightarrow}
	S^q \overset{E_q}\to{\longrightarrow}
	U_1\widehat{K}_2(A)
$$
coincides with $\exp_p$.
By (1), the lifting of $adT/T$ in $({\Bbb I}^{[2]}/{\Bbb J}^{[2]})^1=I\otimes\widehat{\Omega}_B^1/J\otimes\widehat{\Omega}_B^1$ is
$p\widetilde{a}\otimes dT/T$, where $\widetilde{a}$ is a lifting of $a$ to $D$.
Chasing the connecting homomorphism $\delta$,
$$
\minCDarrowwidth{2em}
\CD
	\scriptstyle{0} @>>> \scriptstyle{(J\otimes\widehat{\Omega}_B^1) \oplus D} @>>> \scriptstyle{(I\otimes \widehat{\Omega}_B^1) \oplus D} 
		@>>> \scriptstyle{(I\otimes \widehat{\Omega}_B^1)/(J\otimes \widehat{\Omega}_B^1)} @>>> \scriptstyle{0} \\
	@. @V d VV @V d VV @V d VV \\
	\scriptstyle{0} @>>> \scriptstyle{(D\otimes\widehat{\Omega}_B^2) \oplus (D\otimes\widehat{\Omega}_B^1)} @>>> 
		\scriptstyle{(D\otimes\widehat{\Omega}_B^2) \oplus (D\otimes\widehat{\Omega}_B^1)} 
		@>>> \scriptstyle{0} @>>> \scriptstyle{0} \\
	@. @V d VV @V d VV @V d VV \\
\endCD
$$
(the left column is $\scS(2)$, the middle is $\scS'(2)$ and the right is ${\Bbb I}^{[2]}/{\Bbb J}^{[2]}$);
$p\widetilde{a}dT/T$ in the upper right goes 
to $(pd\widetilde{a}\wedge dT/T,(1-f_2)(p\widetilde{a}\otimes dT/T))$ in the lower left.
By $E_2$, this element goes 
$$
\aligned 
	&E_2\bigl((1-f_2)\bigl(p\widetilde{a}\otimes \frac{d}{T}\bigr)\bigr) 
	=E_2\bigl((1-f_1)(p\widetilde{a})\otimes \frac{dT}{T}\bigr) \\
	&=\{E_1((1-f_1)(p\widetilde{a})),T\} 
	=\bigl\{\exp \circ \bigl(\sum_{n\ge 0}f_1^n\bigr)\circ (1-f_1)(p\widetilde{a}),T \bigr\} \\
	&=\{\exp(pa),T\}.
\endaligned 
$$
in $U_1\widehat{K}_2(A)$.
This is none other than the map $\exp_p$.
\qed\endpf

By Theorem 2 we can calculate the kernel of $\exp_p$.
On the other hand, 
even though $\exp_p$ is not surjective,
the image of $\exp_p$ includes $U_{e+1}\widehat{K}_q(A)$
and we already know $\gr_i\widehat{K}_q(K)$ for $0 \le i \le ep/(p-1)$.
Thus it is enough to calculate the kernel of $\exp_p$ in order to know all $\gr_i\widehat{K}_q(K)$.
Note that to know $\gr_i\widehat{K}_q(K)$,
we may assume that $\zeta_p \in K$, and hence $\widehat{K}_q(A)=U_0\widehat{K}_q(K)$.

\HH 15.2. Computation of the kernel of the exponential homomorphism

\HHH 15.2.1. Modified syntomic complex

We introduce a modification of ${\scS}'(q)$ and calculate it instead of ${\scS}'(q)$.
Let ${\Bbb S}_q$ be the mapping fiber complex of
$$
	1-f_q \: ({\Bbb J}^{[q]})^{\ge q-2} \longrightarrow {\Bbb D}^{\ge q-2}.
$$
Here, for a complex $C^{\cdot}$, we put
$$
	C^{\ge n}=(0 \longrightarrow \dots \longrightarrow 0 \longrightarrow C^n 
		\longrightarrow C^{n+1} \longrightarrow \cdots).
$$ 

By  definition, we have a natural surjection $H^{q-1}({\Bbb S}_q) \rightarrow H^{q-1}(\scS'(q))$,
hence $\psi(H^{q-1}({\Bbb S}_q))=\psi(H^{q-1}(\scS'(q)))$, which is the kernel of $\exp_p$.
\par

To calculate $H^{q-1}({\Bbb S}_q)$, we introduce an $X$-filtration.
Let $0 \le  r \le  2$ and $s=q-r$.
Recall that $B=A_0[[X]]$.
For $i \ge 0$, let $\fil_i (I^{[r]}\otimes_B \widehat{\Omega}_B^s)$ be 
the subgroup of $I^{[r]}\otimes_B \widehat{\Omega}_B^{s}$
generated by the elements 
$$
	\aligned 
	&\biggl\{ X^n  \frac{(X^e)^j}{j!}\frac{p^l}{l!} a\omega  : 
		n+ej \ge i, n\ge 0, j+l \ge r,a \in D, \omega \in \widehat{\Omega}_B^{s} \biggr\} \\ 
	&\cup 
	\biggl\{ X^{n} \frac{(X^e)^j}{j!}\frac{p^l}{l!} a\upsilon \wedge \frac{dX}{X}  : 
		n+ej \ge i, n \ge 1, j+l \ge r, a \in D, \upsilon 
		\in \widehat{\Omega}_B^{s-1} \biggr\}.
	\endaligned 
$$
The map $1-f_q\colon I^{[r]}\otimes \widehat{\Omega}_B^s \rightarrow D \otimes \widehat{\Omega}_B^s$
preserves the filtrations. 
By using the latter we get the following 

\th Proposition 3 

$H^{q-1}(\fil_i{\Bbb S}_q)_i$ form a finite decreasing filtration 
of $H^{q-1}({\Bbb S}_q)$.
Denote
$$
\aligned
&\fil_iH^{q-1}({\Bbb S}_q)=H^{q-1}(\fil_i{\Bbb S}_q), \\
&\gr_iH^{q-1}({\Bbb S}_q)=\fil_iH^{q-1}({\Bbb S}_q)/\fil_{i+1}H^{q-1}
({\Bbb S}_q).
\endaligned 
$$
Then $\gr_i H^{q-1}({\Bbb S}_q)$
$$= 
	\cases
		0 
&\quad (\text{ if }i > 2e) \\
		X^{2e-1}dX \wedge \left(\widehat{\Omega}_{A_0}^{q-3}/p\right) &\quad (\text{ if }i =2e) \\
X^i\left(\widehat{\Omega}_{A_0}^{q-2}/p\right) \oplus X^{i-1}dX\wedge(\widehat{\Omega}_{A_0}^{q-3}/p) 
&\quad (\text{ if } e<i<2e) \\
		X^e\left(\widehat{\Omega}_{A_0}^{q-2}/p\right) 
			\oplus X^{e-1}dX\wedge \left({\frak Z}_1 \widehat{\Omega}_{A_0}^{q-3}\big/ p^2\widehat{\Omega}_{A_0}^{q-3}\right) 
			&\quad (\text{ if } i=e, p\mid e) \\
		X^{e-1}dX\wedge \left({\frak Z}_1 \widehat{\Omega}_{A_0}^{q-3} \big/ p^2\widehat{\Omega}_{A_0}^{q-3}\right)
			&\quad (\text{ if } i=e, p\nmid e) \\
		\left(X^i  
\frac
{\left( p^{\max(\eta_i'-v_p(i),0)}\widehat{\Omega}_{A_0}^{q-2} \cap 
			{\frak Z}_{\eta_i}\widehat{\Omega}_{A_0}^{q-2} \right)+p^2\widehat{\Omega}_{A_0}^{q-2}}
			{p^2\widehat{\Omega}_{A_0}^{q-2}}  \right) \\
		 \qquad \oplus \left( X^{i-1}dX\wedge 
			\frac
{{\frak Z}_{\eta_i}\widehat{\Omega}_{A_0}^{q-3}+p^2\widehat{\Omega}_{A_0}^{q-3} }
{p^2\widehat{\Omega}_{A_0}^{q-3}}
 \right)
	& \quad (\text{ if } 1 \le i <e) \\
		0 &\quad (\text{ if } i=0).
\endcases
$$
Here $\eta_i$ and $\eta_i'$ are the integers which satisfy 
$p^{\eta_i-1}i < e \le p^{\eta_i}i$ and $p^{\eta_i'-1}i-1 < e \le p^{\eta_i'}i-1$ for each $i$,
$$
	{\frak Z}_n\widehat{\Omega}_{A_0}^q = \kr \left( 
		\widehat{\Omega}_{A_0}^q \overset{d}\to{\longrightarrow} \widehat{\Omega}_{A_0}^{q+1}/p^n \right)
$$
for positive $n$,
and ${\frak Z}_n\widehat{\Omega}_{A_0}^q=\widehat{\Omega}_{A_0}^q$ for $n\le 0$.
\endth

\pf Outline of the proof

From the definition of the filtration we have the exact sequence of complexes:
$$
	0 \longrightarrow \fil_{i+1}{\Bbb S}_q \longrightarrow \fil_i{\Bbb S}_q
		\longrightarrow \gr_i{\Bbb S}_q \longrightarrow 0
$$
and this sequence induce a long exact sequence
$$
	\dots \rightarrow H^{q-2}(\gr_i{\Bbb S}_q)
		\rightarrow H^{q-1}(\fil_{i+1}{\Bbb S}_q)
		\rightarrow H^{q-1}(\fil_i{\Bbb S}_q)
		\rightarrow H^{q-1}(\gr_i{\Bbb S}_q)
		\rightarrow \cdots.
$$
The group $H^{q-2}(\gr_i{\Bbb S}_q)$ is
$$
	H^{q-2}(\gr_i{\Bbb S}_q)
	=\kr
	\binom{\gr_iI^{[2]}\otimes \widehat{\Omega}_B^{q-2} \longrightarrow
			(\gr_iI\otimes \widehat{\Omega}_B^{q-1})\oplus (\gr_iD\otimes\widehat{\Omega}_B^{q-2})}
		{\! \! \! \! \! \! \! \! \! \! x \longmapsto (dx, (1-f_q)x )}.
$$
The map $1-f_q$ is equal to $1$ if $i\ge 1$ and $1-f_q \: p^2\widehat{\Omega}_{A_0}^{q-2}\rightarrow \widehat{\Omega}_{A_0}^{q-2}$
if $i=0$,
thus they are all injective.
Hence $H^{q-2}(\gr_i{\Bbb S}_q)=0$ for all $i$
and we deduce that
$H^{q-1}(\fil_i{\Bbb S}_q)_i$ form a decreasing filtration on $H^{q-1}({\Bbb S}_q)$.
\par

Next, we have to calculate $H^{q-2}(\gr_i{\Bbb S}_q)$.
The calculation is easy but there are many cases which depend on $i$,
so we omit them.
For more detail, see \cite{N2}.

Finally, we have to compute the image of the last arrow of the exact sequence
$$
	0
		\longrightarrow H^{q-1}(\fil_{i+1}{\Bbb S}_q)
		\longrightarrow H^{q-1}(\fil_i{\Bbb S}_q)
		\longrightarrow H^{q-1}(\gr_i{\Bbb S}_q)
$$
because it is not surjective in general.
Write down the complex $\gr_i{\Bbb S}_q$:
$$
	\dots
	\rightarrow (\gr_i I\otimes \widehat{\Omega}_B^{q-1})\oplus (\gr_i D\otimes\widehat{\Omega}_B^{q-2})
	\overset{d}\to{\rightarrow} (\gr_i D\otimes \widehat{\Omega}_B^{q})\oplus (\gr_i D\otimes\widehat{\Omega}_B^{q-1})
	\rightarrow \cdots,
$$
where the first term is the degree $q-1$ part and the second term is the degree $q$ part.
An element $(x,y)$ in the first term which is mapped to zero by $d$
comes from $H^{q-1}(\fil_i{\Bbb S}_q)$ if and only if
there exists $z \in \fil_i D\otimes\widehat{\Omega}_B^{q-2}$ such that 
$z\equiv y$ modulo $\fil_{i+1} D\otimes\widehat{\Omega}_B^{q-2}$ and
$$
	\sum_{n\ge 0}f_q^n(dz) \in \fil_i I\otimes \widehat{\Omega}_B^{q-1}.
$$
From here one deduces Proposition 3.
\qed\endpf

\HHH 15.2.2. Differential modules

Take a prime element $\pi$ of $K$ such that   $\pi^{e-1}d\pi=0$. 
We assume that $p \nmid e$ in this subsection.
Then we have
$$\split
	&\widehat{\Omega}_A^q \simeq \biggl( \bigoplus_{i_1<i_2<\dots<i_q} 
		A\frac{dT_{i_1}}{T_{i_1}}\wedge \dots\wedge \frac{dT_{i_q}}{T_{i_q}} \biggr) \\
	&\quad \oplus \biggl( \bigoplus_{i_1<i_2<\dots<i_{q-1}} 
		A/(\pi^{e-1})\frac{dT_{i_1}}{T_{i_1}}\wedge \dots\wedge \frac{dT_{i_{q-1}}}{T_{i_{q-1}}}\wedge d\pi \biggr),
\endsplit
$$
where $\{T_i\}={\Bbb T}$.
We introduce a filtration on $\widehat{\Omega}_A^q$ as
$$
	\fil_i \widehat{\Omega}_A^q =
		\cases
			\widehat{\Omega}_A^q &\quad (\text{ if }i=0) \\
			\pi^i \widehat{\Omega}_A^q +\pi^{i-1}d\pi \wedge \widehat{\Omega}_A^{q-1} &\quad (\text{ if }i \ge 1).
		\endcases
$$
The subquotients are
$$
\split
	&\gr_i\widehat{\Omega}_A^q = \fil_i\widehat{\Omega}_A^q/\fil_{i+1}\widehat{\Omega}_A^q \\
		&=\cases
			\Omega_F^q &\quad (\text{ if } i=0 \text{ or } i\ge e) \\
			\Omega_F^q \oplus \Omega_F^{q-1} &\quad (\text{ if } 1 \le i < e),
		\endcases
\endsplit
$$
where the map is
$$
\aligned 
	&\Omega_F^q  \ni \omega \longmapsto \pi^i\widetilde{\omega} \in  \pi^i\widehat{\Omega}_A^q\\
	&\Omega_F^{q-1} \ni \omega \longmapsto \pi^{i-1}d\pi \wedge \widetilde{\omega} \in \pi^{i-1}d\pi \wedge \widehat{\Omega}_A^{q-1}.
\endaligned 
$$
Here $\widetilde{\omega}$ is the lifting of $\omega$.
Let $\fil_i(\widehat{\Omega}_A^q/pd\widehat{\Omega}_A^{q-1})$ be the image of $\fil_i\widehat{\Omega}_A^q$ in 

\noindent $\widehat{\Omega}_A^q/pd\widehat{\Omega}_A^{q-1}$.
Then we have the following:

\th Proposition 4

For $j \ge 0$, 
$$
	\gr_j \left(\widehat{\Omega}_A^q/pd\widehat{\Omega}_A^{q-1} \right) =
	\cases
		\Omega_F^q \quad &(j=0) \\
		\Omega_F^q \oplus \Omega_F^{q-1} \quad &(1 \le j <e) \\
		\Omega_F^q/{B}_l^q \quad &(e \le j),
	\endcases
$$
where $l$ be the maximal integer which satisfies $j-le \ge 0$.
\endth

\pf Proof

If $1 \le j <e$, $\gr_j\widehat{\Omega}_A^q=\gr_j(\widehat{\Omega}_A^q/pd\widehat{\Omega}_A^{q-1})$ because $pd\widehat{\Omega}_A^{q-1} \subset \fil_e\widehat{\Omega}_A^q$.
Assume that $j \ge e$ and let $l$ be as above.
Since $\pi^{e-1}d\pi=0$, $\widehat{\Omega}_A^{q-1}$ is generated by 
elements $p\pi^id\omega$ for $0 \le i <e$ and $\omega \in \widehat{\Omega}_{A_0}^{q-1}$.
By \cite{I} (Cor.\ 2.3.14), $p\pi^i d\omega \in \fil_{e(1+n)+i}\widehat{\Omega}_A^q$ if and only if
the residue class of $p^{-n}d\omega$ belongs to $B_{n+1}$.
Thus $\gr_j(\widehat{\Omega}_A^q/pd\widehat{\Omega}_A^{q-1}) \simeq \Omega_F^q/{B}_l^q$.
\qed\endpf 

By definition of the filtrations, 
$\exp_p$ preserves the filtrations on $\widehat{\Omega}_A^{q-1}/pd\widehat{\Omega}_A^{q-2}$ and $\widehat{K}_q(K)$.
Furthermore,
$\exp_p \colon \gr_i (\widehat{\Omega}_A^{q-1}/pd\widehat{\Omega}_A^{q-2}) \rightarrow \gr_{i+e}K_q(K)$
is surjective and its kernel is the image of 
$\psi(H^{q-1}({\Bbb S}_q)) \cap \fil_i(\widehat{\Omega} _A^{q-1}/pd\widehat{\Omega}_A^{q-2})$ 
in $\gr_i (\widehat{\Omega}_A^{q-1}/pd\widehat{\Omega}_A^{q-2})$.
Now we know both $\widehat{\Omega} _A^{q-1}/pd\widehat{\Omega}_A^{q-2}$ and $H^{q-1}({\Bbb S}_q)$ explicitly,
thus we shall get the structure of $K_q(K)$ by calculating $\psi$. 
But $\psi$ does not preserve the filtration of $H^{q-1}({\Bbb S}_q)$,
so it is not easy to compute it.
For more details, see \cite{N2}, especially sections 4-8 of that paper.
After completing these calculations, we get the result in (vi) in the introduction.

\rk Remark

Note that if $p\mid e$, the structure of $\widehat{\Omega}_A^{q-1}/pd\widehat{\Omega}_A^{q-2}$ 
is much more complicated.
For example, if $e=p(p-1)$, and if $\pi^e=p$, then $p\pi^{e-1}d\pi=0$.
This means the torsion part of $\widehat{\Omega}_A^{q-1}$ is larger than in the the case where $p\nmid e$.
Furthermore, if $\pi^{p(p-1)}=pT$ for some $T \in {\Bbb T}$, 
then $p\pi^{e-1}d\pi=pdT$, this means that $d\pi$ is not a torsion element.
This complexity makes it difficult to describe the structure of $K_q(K)$
in the case where $p \mid e$.
\endrk

\HHH Appendix. The mapping fiber complex

{}

This subsection is only a note on homological algebra  
to introduce the mapping fiber complex.
The mapping fiber complex is the degree $-1$ shift of the mapping cone complex.
\par

Let $C^{\cdot} \overset{f}\to{\rightarrow} D^{\cdot}$ be a morphism of non-negative cochain complexes.
We denote the degree $i$ term of $C^{\cdot}$ by $C^i$.

Then the mapping fiber complex $\operatorname{MF}(f)^{\cdot}$ is defined as follows.
$$
\split
	&\operatorname{MF}(f)^i = C^i \oplus D^{i-1}, \\
	&\text{differential} \quad d \colon~C^i \oplus D^{i-1} \longrightarrow C^{i+1} \oplus D^i \\
	& \kern 9em ~(x,y) \longmapsto (dx,f(x)-dy).
\endsplit
$$
By definition,
we get an exact sequence of complexes:
$$
	0 \longrightarrow D[-1]^{\cdot} \longrightarrow \operatorname{MF}(f)^{\cdot} 
	\longrightarrow C^{\cdot} \longrightarrow 0,
$$
where $D[-1]^{\cdot}=(0 \rightarrow D^0 \rightarrow D^1 \rightarrow \cdots)$ 
	(degree $-1$ shift of $D^{\cdot}$.)
\par

Taking cohomology, we get a long exact sequence
$$
	\cdots \rightarrow H^i(\operatorname{MF}(f)^{\cdot}) \rightarrow H^i(C^{\cdot}) \rightarrow H^{i+1}(D^{\cdot}[-1])
	\rightarrow H^{i+1}(\operatorname{MF}(f)^{\cdot}) \rightarrow \cdots,
$$
which is the same as the following exact sequence
$$
	\cdots \rightarrow H^i(\operatorname{MF}(f)^{\cdot}) \rightarrow H^i(C^{\cdot}) 
	\overset{f}\to{\rightarrow} H^{i}(D^{\cdot})
	\rightarrow H^{i+1}(\operatorname{MF}(f)^{\cdot}) \rightarrow \cdots .
$$
\par

\Bib References

\rf{B}
S. Bloch, 
 Algebraic {$K$}-theory and crystalline cohomology, 
 {Publ. Math. IHES},  47(1977), 187--268.

\rf{BK}
S. Bloch and K. Kato,
 $p$-adic \'etale cohomology, 
 {Publ. Math. IHES}, 63(1986), 107--152.

\rf{BT}
H. Bass and J. Tate, 
 The {Milnor} ring of a global field, 
 In {Algebraic $K$-theory II}, Lect. Notes Math. 342,
Springer-Verlag, Berlin, 1973, 349--446.

\rf{F}
I. Fesenko, 
Abelian local $p$-class field theory,
Math. Ann., 301 (1995), pp. 561--586.

\rf{G}
J. Graham, 
 Continuous symbols on fields of formal power series, 
 In {Algebraic $K$-theory II}, Lect. Notes Math. 342,
Springer-Verlag, Berlin, 1973, 474--486.

\rf{I}
L. Illusie, 
 Complexe de {de Rham--Witt} et cohomologie cristalline,
 {Ann. Sci. Ecole Norm. Sup. } 12(1979), 501--661.

\rf{Kh}
B. Kahn, 
 L'anneau de {Milnor} d'un corps local \`a corps residuel parfait.
 {Ann. Inst. Fourier 34}, 4(1984), 19--65. 

\rf{Kt1} 
K. Kato, 
 A generalization of local class field theory by using {$K$}-groups. 
  {II}, 
 {J. Fac. Sci. Univ. Tokyo} 27(1980), 603--683.

\rf{Kt2}
K. Kato, 
 On $p$-adic vanishing cycles (applications of ideas of 
  {Fontaine--Messing}). 
 {Adv. Stud. in Pure Math.} 10(1987), 207--251.

\rf{Ku1} 
M. Kurihara, 
 Abelian extensions of an absolutely unramified local field with
  general residue field, 
 {Invent. Math.} 93(1988), 451--480. 

\rf{Ku2} 
M. Kurihara, 
 The {Milnor $K$-groups} of a local ring over a ring of the $p$-adic
  integers, 
 to appear in the volume on Ramification theory for arithmetic
  schemes. Lminy 1999, ed.~B. Erez.

\rf{Ku3} 
M. Kurihara, 
 On the structure of the {Milnor $K$-group} of a certain complete
  discrete valuation field, 
 preprint.

\rf{Ku4} 
M. Kurihara, 
A note on $p$-adic etale cohomology,
{Proc. Japan Acad. (A) } 63 (1987), 275--278.

\rf{N1}
J. Nakamura, 
 On the structure of the {Milnor $K$-groups} of some complete discrete
  valuation fields, 
{$K$-Theory 19} (2000), 269--309.

\rf{N2}
J. Nakamura, 
 On the {Milnor $K$-groups} of complete discrete valuation fields,  
{Doc. Math.} 5(2000), 151--200.

\rf{P}
A. N. Parshin,
 Local class field theory,
 Trudy Mat. Inst. Steklov  (1984);
English translation in Proc. Steklov Inst. Math. 165 (1985),  no. 3, pp. 157--185.

\rf{Z}
I. Zhukov, 
 Milnor and topological {$K$}-groups of multidimensional complete
  fields, Algebra i Analiz  (1997); 
English  translation in 
 St. Petersburg Math. J.  9(1998), 69--105.
\endBib

 \Coordinates

Department of Mathematics \ 
University of Tokyo

3-8-1 Komaba Meguro-Ku Tokyo 153-8914 Japan

E-mail: jinya\@ms357.ms.u-tokyo.ac.jp
\endCoordinates

\vfill
\pagebreak

\bye